\documentclass{amsart}
\usepackage{amssymb, amsmath, color}

\def\blue{\color{blue}}

\newtheorem{theorem}{Theorem}

\newtheorem{remark}[theorem]{Remark}

\begin{document}

\title{Three-dimensional Lorentzian homogeneous  Ricci solitons}
\author{M. Brozos-V\'{a}zquez $\,$ G. Calvaruso $\,$ E. Garc\'{\i}a-R\'{\i}o $\,$ S. Gavino-Fern\'{a}ndez}
\address{MBV: E. U. Polit\'ecnica, Department of Mathematics, University of A Coru\~na,
Ferrol,
Spain} \email{mbrozos@udc.es} {\blue \address{GC: Dipartimento di
Matematica "E. De Giorgi", Universit\`{a} del Salento, Lecce, Italy}
\email{giovanni.calvaruso@unisalento.it}}
\address{EGR-SGF: Faculty of Mathematics,
University of Santiago de Compostela, 15782 Santiago de Compostela,
Spain} \email{eduardo.garcia.rio@usc.es $\,\,$ sandra.gavino@usc.es}
\thanks{Supported by projects MTM2006-01432 and PGIDIT06PXIB207054PR (Spain) {and MIUR, PRIN 2007 (Italy)}. }

\begin{abstract}
We study three-dimensional Lorentzian homogeneous Ricci solitons, proving the existence
of shrinking, expanding and steady Ricci solitons.
For all the   non-trivial examples, the Ricci operator is not diagonalizable and has three equal eigenvalues.
\end{abstract}

\maketitle

\section{Introduction}

A \emph{Ricci soliton} is a pseudo-Riemannian manifold $(M,g)$ which admits
a smooth vector field $X$ on $M$ such that
\begin{equation}\label{soliton}
\mathcal{L}_Xg+Ric=\lambda g,
\end{equation}
where $\mathcal{L}_X$ denotes the Lie derivative in the direction of
$X,$ $Ric$ denotes the Ricci tensor and $\lambda$ is a real number.
A Ricci solition is said to be a \emph{shrinking}, \emph{steady} or
\emph{expanding}, respectively, if $\lambda>0$, $\lambda=0$ or $\lambda<0${\blue .}

{ The description of} Ricci solitons {can be regarded as a first step} in understanding the
\emph{Ricci flow}, since they are the fixed points of the flow. Moreover, they are
important in understanding singularities of the Ricci flow.
Under suitable conditions, type I singularity models correspond to shrinking solitons,
type II models correspond to steady Ricci solitons while type III models correspond to expanding
Ricci solitons. We {refer} to \cite{Cao} for a survey and further references on the geometry of Ricci solitons.

Theoretical physicists have also been looking into the equation of Ricci solitons in
relation with String Theory. {A seminal} contribution in this direction is due to
Friedan \cite{Fr} (see also \cite{AW, Kh} for a discussion of physical aspects of the Ricci flow).
{Although it was first} introduced and studied in { a} Riemannian { context}, the Ricci soliton equation \eqref{soliton} {
is  currently being investigated} in pseudo-Riemannian settings, with special attention to the
Lorentzian case \cite{Cas,PT}. The Ricci soliton equation may also give some insight into the
general study of {\em Einstein field equations}, of which \eqref{soliton} is a special case.

As proved in \cite{Cerbo} (see also \cite{H, Pa}),
three-dimensional Lie groups do not admit left-invariant Riemannian Ricci solitons.
In this paper, we study the corresponding existence problem in Lorentzian { signature}.
We shall conclude that the Lorentzian case is much richer,
allowing the existence of expanding, steady and shrinking left-invariant Ricci solitons.

Three-dimensional locally homogeneous Lorentzian manifolds are either locally symmetric or locally isometric to a
three-dimensional Lie group equipped with a left-invariant Lorentzian
metric. Moreover, three-dimensional locally symmetric Lorentzian manifolds which are not of constant
sectional curvature, are either locally isometric to a Lorentzian product of a real line and a surface of
constant Gauss curvature, or they are Walker manifolds with a two-step nilpotent Ricci operator~\cite{calvaruso}.

It is clear that, in addition to Einstein spaces, products $N^k(c)\times\mathbb{R}$ with
$N^k(c)$ of constant sectional curvature are Ricci solitons, both in the Riemannian { and the Lorentzian case}.
In fact, it suffices to consider the
Gaussian soliton on the $\mathbb{R}$-factor, scaling it by the constant value of the sectional
curvature. An immediate calculation shows that the gradient of $f(t)=\frac{1}{2}\lambda t^2$ defines a
Ricci soliton on $N^k(c)\times\mathbb{R}$ for $\lambda=(k-1)c$.
Ricci solitons on Walker manifolds with nilpotent Ricci operator are considered in Section~\ref{s:walker},
proving the existence of expanding, steady and shrinking  locally symmetric Ricci solitons.

The existence of left-invariant Lorentzian  Ricci solitons on Lie groups is analyzed in Section~\ref{s:lie}.
The following result summarizes the classification of { Lorentzian} homogeneous Ricci solitons in dimension three:

\begin{theorem}\label{th:1}
Let $M$ be a three-dimensional simply connected, complete  homogeneous manifold.
$M$ is a Ricci soliton if and only if $M$ is one of the spaces listed below:
\begin{enumerate}
\item[i)] a space of constant curvature $N^3(c)$;
\item[ii)] a Lorentzian product $\mathbb{R}\times N^2(c)$;
\item[iii)] a symmetric Walker manifold;
\item[iv)] a unimodular Lie group $G$ with one of the following Lie algebras:
\begin{enumerate}
\item[]
$$
\begin{array}{lll}
(iv.1) & [e_1,e_2]=\frac{1}{2}e_2-(\beta-\frac{1}{2}) e_3,
& [e_1,e_3]=-(\beta+\frac{1}{2}) e_2-\frac{1}{2} e_3, {\vphantom{\displaystyle{\frac{1}{B}}}}\\
&[e_2,e_3]=\alpha e_1,&
\end{array}
$$
with either $\alpha=0$ or $\alpha=\beta \neq 0$.
{ If $\alpha=0$ then $G=E(1,1)$, while if $\alpha =\beta \neq 0$ then $G=(1,2)$ or $G=SL(2,\mathbb R)$}.
\item[]
$$\begin{array}{ll}
(iv.2) \quad & [e_1,e_2]=-\frac{1}{\sqrt{2}}e_1-\alpha e_3, \qquad
[e_1,e_3]=-\frac{1}{\sqrt{2}} e_1-\alpha e_2, \\ &
{[e_2,e_3]=\alpha e_1+\frac{1}{\sqrt{2}}e_2 -\frac{1}{\sqrt{2}}e_3.}
\end{array}
$$
{ If $\alpha =0$ then $G=E(1,1)$, while if $\alpha \neq 0$ then either $G=(1,2)$ or $G=SL(2,\mathbb R)$}.
\end{enumerate}
\item[v)] a non-unimodular Lie group $G$ with Lie algebra given by
\[
\begin{array}{rcl}
[e_1,e_2]\!&\!=\!&\!-\frac{1}{\sqrt{2}}\left(\alpha e_1+\frac{1}{\sqrt{2}}\beta( e_2+ e_3)\right),
\quad
[e_1,e_3]\!\!=\!\!\frac{1}{\sqrt{2}}\left(\alpha e_1+\frac{1}{\sqrt{2}}\beta( e_2+ e_3)\right),\\
\noalign{\smallskip}
[e_2,e_3]\!&\!=\!&\!\frac{1}{\sqrt{2}}\delta(e_2 +e_3).
\end{array}
\]
\end{enumerate}
In { cases iv) and v)}, $\{ e_1, e_2, e_3\}$ is an orthonormal basis of the corresponding Lie algebra, of signature $(++-)$.
\end{theorem}

From now on, by a \emph{non-trivial} Ricci soliton we shall mean { a Ricci soliton which} is neither
Einstein nor a product $\mathbb{R}\times N^k(c)$.
The explicit description of three-dimensional non-trivial homogeneous
Lorentzian Ricci solitons { is given in Sections~\ref{s:lie} and \ref{s:walker},}
where the causal character of vector fields defining these Ricci solitons { is also} discussed.

We briefly recall that the Ricci operator $\hat{Ric}$, being self-adjoint,
is always diagonalizable in the Riemannian case, while at each point of a Lorentzian manifold
four different cases can occur, known as \emph{Segre types}.
In dimension three, the possible cases are the following \cite{HMP}:
\begin{enumerate}
\item[(Ia)] \emph{Segre type $\{11,1\}$}:
$\hat{Ric}$ is symmetric and { hence} diagonalizable.
The comma separates the spacelike and timelike eigenvectors.
In the degenerate case, at least two of the Ricci eigenvalues coincide.
\item[(Ib)] \emph{Segre type $\{1 z \bar{z} \}$}: $\hat{Ric}$ has one real and two complex conjugate eigenvalues.
\item[(II)] \emph{Segre type $\{21 \}$}:
$\hat{Ric}$ has two real eigenvalues ({ which coincide} in the degenerate case),
one of which has multiplicity two and each associated to a one-dimensional eigenspace.
\item[(III)] \emph{Segre type $\{3\}$}:
$\hat{Ric}$ has three equal eigenvalues, associated to a one-dimen\-sional eigenspace.
\end{enumerate}
Segre types of the Ricci operator of a three-dimensional homogeneous Lorentzian manifold were discussed
in \cite{CK} (see also the description included { in} Sections~\ref{s:lie} and \ref{s:walker}). Taking into account the classification given in Theorem  \ref{th:1}, we have at once the following

\begin{theorem}\label{th:2}
A complete and simply connected three-dimensional homogeneous Lorentzian manifold is a non-trivial Ricci soliton if and only
if the Ricci operator
$\hat{Ric}$ is not diagonalizable and has exactly three equal eigenvalues, that is,
$\hat{Ric}$ is either of Segre type $\{3\}$ or of degenerate Segre type $\{21\}$.
\end{theorem}

\section{Three-dimensional Homogeneous Lorentzian Manifolds}\label{s:lie}

A connected, complete and simply connected
three-dimensional  homogeneous Lorentzian manifold is a Lie group \cite{calvaruso}.
For the sake of completeness we include a brief description of
three-dimensional unimodular and non-unimodular Lie groups. Theorems
\ref{th:1} and \ref{th:2} will follow from the subsequent analysis.
{ Theorem
\ref{th:1} iv) and v) and Theorem~\ref{th:2}} will follow from the subsequent analysis.

\subsection{Unimodular Lie groups}

Let $\times$ denote the Lorentzian vector product on
$\mathbb{R}^3_1$ induced by the product of the para-quaternions
(i.e., $e_1\times e_2=-e_3$, $e_2\times e_3=e_1$, $e_3\times
e_1=e_2$, where $\{ e_1,e_2,e_3\}$ is an orthonormal basis of
signature $(++-)$). The Lie bracket $[,]$ defines the corresponding
Lie algebra $\mathfrak{g}$, which is unimodular if and only if the
endomorphism $L$ defined by $[Z,Y]=L(Z\times Y)$ is self-adjoint
\cite{Rahmani}. Considering the different Segre types of
$L$, we have the following four classes of unimodular
three-dimensional Lie algebras (we follow notation in
\cite{GR-HB-VL}):

\subsubsection{\underline{\emph{Segre type $\{11,1\}$}}.}

If $L$ is diagonalizable with eigenvalues $\{\alpha,\beta,\gamma\}$
with respect to an orthonormal basis $\{e_1,e_2,e_3\}$ of signature
$(++-)$, the corresponding Lie algebra is given by
\begin{equation}\label{Ia}
(\mathfrak{g}_{Ia}): \qquad\qquad [e_1,e_2]=-\gamma e_3, \quad
[e_1,e_3]=-\beta e_2, \quad [e_2,e_3]=\alpha e_1.
\end{equation}
Up to symmetries, the only nonvanishing components of the curvature tensor are given by
\[\begin{array}{l}
R_{1221}=\frac{1}{4} \left(\alpha ^2+\beta ^2-3 \gamma ^2-2
\alpha\beta  +2 \alpha\gamma
 +2 \beta\gamma \right), \\[0.05in]
R_{1313}=\frac{1}{4} \left(\alpha ^2-3 \beta ^2+\gamma ^2+2
\alpha\beta   -2 \alpha\gamma
 +2 \beta  \gamma \right), \\[0.05in]
R_{2332}=\frac{1}{4} \left(3 \alpha ^2-\beta^2-\gamma^2-2
\alpha\beta -2\alpha\gamma +2\beta\gamma\right),
\end{array}\]
and the Ricci operator is diagonalizable (that is, of Segre type $\{11,1\}$) with respect to the basis
$\{e_1,e_2,e_3\}$ with eigenvalues
\begin{equation}\label{eq:eig-ricci-Ia}\textstyle
\lambda_1=\frac{1}{2}((\beta-\gamma)^2-\alpha^2),\quad
\lambda_2=\frac{1}{2}((\alpha-\gamma)^2-\beta^2),\quad
\lambda_3=\frac{1}{2}((\alpha-\beta)^2-\gamma^2)\,.
\end{equation}

For an arbitrary vector $X=\sum X_i e_i$, from equation (\ref{Ia}) we get
\[
(\mathcal{L}_Xg)= \left(
\begin{array}{ccc}
0 & X_3 (\alpha-\beta) & X_2 (\gamma-\alpha)\\
X_3 (\alpha-\beta) & 0 & X_1 (\beta-\gamma)\\
X_2 (\gamma-\alpha) & X_1 (\beta-\gamma) & 0
\end{array}\right).
\]
Hence, by \eqref{soliton},  there exist a Ricci soliton of this type if and only if the following system of equations is satisfied:
\begin{equation}\label{Iaeq}
\left\{
\begin{array}{l}
(\beta-\gamma)^2-\alpha^2=2\lambda{,}
\\
\noalign{\smallskip}
(\alpha-\gamma)^2-\beta^2=2\lambda{,}
\\
\noalign{\smallskip}
(\alpha-\beta)^2-\gamma^2=2\lambda{,}\\
\noalign{\smallskip}
X_1(\beta-\gamma)=0{,}
\\
\noalign{\smallskip}
X_2(\alpha-\gamma)=0{,}
\\
\noalign{\smallskip}
X_3(\alpha-\beta)=0{.}
\end{array}
\right.
\end{equation}

Now, from (\ref{eq:eig-ricci-Ia}), it is clear that any solution of
(\ref{Iaeq}) gives rise to an Einstein metric. Therefore \emph{there
are no homogeneous non-trivial Ricci solitons of Segre type $\{11,1\}$}.

\subsubsection{\underline{\emph{Segre type $\{1 z \bar{z} \}$}}.}

Assume $L$ has a complex eigenvalue. Then
\[
L=\left(\begin{array}{ccc}
\alpha &0&0\\
0&\gamma&-\beta\\
0 &\beta &\gamma
\end{array}\right), \qquad \beta\neq 0,
\]
with respect to an orthonormal basis $\{e_1,e_2,e_3\}$ of signature $(++-)$.
The corresponding Lie algebra is given by
\[
(\mathfrak{g}_{Ib}): \quad\quad [e_1,e_2]=\beta e_2-\gamma e_3,
\quad [e_1,e_3]=-\gamma e_2-\beta e_3, \quad [e_2,e_3]=\alpha e_1.
\]

The non-zero components of the curvature tensor (up to symmetries)
 are
\[\textstyle
R_{1221}=R_{1313}=\frac{1}{4} \left(\alpha ^2+4 \beta ^2\right),
\quad R_{2332}=\textstyle{\frac{3}{4}}\alpha ^2+\beta
^2-\alpha\gamma , \quad R_{1231}=\beta(\alpha-2\gamma).
\]
The Ricci operator $\hat{Ric}$, with respect to the basis $\{e_1,e_2,e_3\}$, is described as follows:
\[
\hat{Ric}=\left(\begin{array}{ccc}
-\frac{1}{2}(\alpha^2+4\beta^2) &0&0\\
0&\frac{1}{2}\alpha(\alpha-2\gamma)&-\beta(\alpha-2\gamma)\\
0 &\beta(\alpha-2\gamma) &\frac{1}{2}\alpha(\alpha-2\gamma)
\end{array}\right), \qquad \beta\neq 0.
\]
Hence, $\hat{Ric}$ is of Segre type  $\{1z \bar z\}$ if $\alpha\neq 2\gamma$ and $\{11,1\}$ if $\alpha=2\gamma$.
For $X=\sum X_i e_i$, one has
\[
(\mathcal{L}_Xg)= \left(
\begin{array}{ccc}
0 & X_2\beta+X_3(\alpha-\gamma) & X_3\beta+X_2(\gamma-\alpha)\\
X_2\beta+X_3(\alpha-\gamma) & -2X_1\beta & 0\\
X_3\beta+X_2(\gamma-\alpha) & 0& -2X_1\beta
\end{array}\right){,}
\]
and thus, we have a homogeneous Ricci soliton of
Segre type $\{1 z \bar{z} \}$ if and only if
\begin{equation}\label{Ibeq}
\left\{
\begin{array}{l}
\alpha^2+4\beta^2=-2\lambda{,}
\\
\noalign{\smallskip}
\alpha^2-2\alpha\gamma-4 X_1\beta=2\lambda {,}\\
\noalign{\smallskip}
\alpha^2-2\alpha\gamma+4X_1 \beta =2\lambda{,}
\\
\noalign{\smallskip}
X_3(\alpha-\gamma)+X_2\beta=0{,}
\\
\noalign{\smallskip}
X_2(\alpha-\gamma)-X_3\beta=0{,}
\\
\noalign{\smallskip}
\beta(\alpha-2\gamma)=0{.}
\end{array}
\right.
\end{equation}

Since $\beta\neq 0$, the last equation in \eqref{Ibeq} gives $\alpha-2\gamma=0$.
Hence, the second and third equations simplify to $-4 X_1\beta=2\lambda$ and
$4 X_1\beta=2\lambda$, respectively, which imply $X_1=\lambda=0$.
Finally, from the first equation one gets that there are no
solutions of (\ref{Ibeq}) with $\beta\neq 0$. Therefore \emph{there
are no homogeneous Ricci solitons of Segre type $\{1 z \bar{z} \}$}.

\subsubsection{\underline{\emph{Segre type $\{21 \}$}}.}

Assume $L$ has a double root of its minimal polynomial. Then, with
respect to an orthonormal basis $\{e_1,e_2,e_3\}$ of signature
$(++-)$, one has
\[
L=\left(\begin{array}{ccc}
\alpha & 0&0\\
     0  &\frac{1}{2}+\beta&-\frac{1}{2}\\
     0  & \frac{1}{2}&-\frac{1}{2}+\beta
\end{array}\right)
\]
and the corresponding Lie algebra is given by
\[
\textstyle (\mathfrak{g}_{II}):\quad
[e_1,e_2]=\frac{1}{2}e_2-(\beta-\frac{1}{2}) e_3, \quad
[e_1,e_3]=-(\beta+\frac{1}{2}) e_2-\frac{1}{2} e_3, \quad
[e_2,e_3]=\alpha e_1.
\]
The non-zero components of the curvature tensor are given by
\[
\begin{array}{ll}
\textstyle R_{1221}=\frac{1}{4} \left(\alpha^2 -2\alpha +4 \beta
\right), & R_{1313}=\frac{1}{4} \left(\alpha^2 +2\alpha
-4\beta\right),
\\[0.05in]
R_{2332}=\frac{1}{4} \alpha (3\alpha -4 \beta), &\quad
R_{1231}=\frac{1}{2}\alpha -\beta.
\end{array}
\]
Hence the Ricci operator takes the form
\begin{equation}\label{RicII}
\hat{Ric}=\left(\begin{array}{ccc}
-\frac{1}{2}\alpha^2 &0&0\\
0&\frac{1}{2}({\alpha+1})(\alpha-2\beta)&-\frac{1}{2}\alpha +\beta\\
0 &\frac{1}{2} \alpha -\beta &\frac{1}{2}(\alpha-1)(\alpha-2\beta)
\end{array}\right),
\end{equation}
with eigenvalues $\lambda_1=-\frac{1}{2}\alpha^2$ and
$\lambda_2=\lambda_3=\frac{1}{2}\alpha(\alpha-2\beta)$.
Thus, $\hat{Ric}$ is of Segre type $\{21\}$, degenerate if and only if $\alpha=0$ or $\alpha=\beta$.

For a vector field $X=\sum X_i e_i$, we get
\[
(\mathcal{L}_Xg)= \left(
\begin{array}{ccc}
0 & a_{12} & a_{13}\\
a_{12} & -X_1 & X_1\\
a_{13} & X_1& -X_1
\end{array}\right).
\]
where $\displaystyle a_{12}=\frac{1}{2}(X_2+X_3(2\alpha-2\beta-1))$ and
$\displaystyle a_{13}=\frac{1}{2}(X_3+X_2(2\beta-2\alpha-1))$.
{Necessary} and sufficient conditions for the existence of a
homogeneous Ricci soliton of Segre type $\{21 \}$ are then given by
\begin{equation}\label{IIeq}
\left\{
\begin{array}{l}
\alpha^2=-2\lambda{,}
\\
\noalign{\smallskip}
\alpha^2-2\alpha\beta+\alpha-2\beta-2X_1=2\lambda{,}\\
\noalign{\smallskip}
\alpha^2-2\alpha\beta-\alpha+2\beta+2X_1=2\lambda{,}
\\
\noalign{\smallskip}
\alpha-2\beta-2X_1=0{,}\\
\noalign{\smallskip}
(2\alpha-2\beta)X_3+X_2-X_3=0{,}
\\
\noalign{\smallskip}
(2\alpha-2\beta)X_2+X_2-X_3=0{.}
\end{array}
\right.
\end{equation}

From the second and forth equation in \eqref{IIeq} one gets $\alpha^2-2\alpha\beta-2\lambda=0$. Replacing into
the first equation, we then obtain $\alpha(\alpha-\beta)=0$.
Hence, either $\alpha=0\neq \beta$ or $\alpha=\beta\neq 0$.
(We excluded the case $\alpha=\beta=0$, since by \eqref{RicII} this corresponds to a flat manifold.)

\underline{First case: $\alpha=0 \neq \beta$.}
From the first equation in (\ref{IIeq}) one gets $\lambda=0$,  the last two equations give
$X_2=X_3$ and the forth equation yields $X_1=-\beta$.
Therefore, the (spacelike) vector field
\begin{equation}\label{XIIa}
X=-\beta\, e_1
\end{equation}
defines a homogeneous (steady) Ricci soliton. By \eqref{RicII},
the Ricci operator is two-step nilpotent but nonvanishing (since $\beta\neq 0$),
that is, of degenerate Segre type $\{21\}$ with eigenvalue equal to zero.

\underline{Second case: $\alpha=\beta\neq 0$.}
In this case, one easily gets  from (\ref{IIeq}) that $\lambda=-\frac{1}{2}\beta^2$,
that $X_1=-\frac{1}{2}\beta$ and that $X_2=X_3$; thus, \emph{there exist a
one-parameter family of homogeneous expanding Ricci solitons}, given by
\begin{equation}\label{XIIb}
X=-\frac{1}{2}\beta \, e_1 + \delta \, e_2 + \delta \, e_3, \quad \delta\in \mathbb{R}.
\end{equation}
Note that the causality of $X$ is again fixed and
one can only find examples of solitons for $X$ spacelike but not null or timelike.
Since $\alpha=\beta \neq 0$, \eqref{RicII} yields that the Ricci operator is of degenerate Segre type $\{21\}$,
with one non-zero eigenvalue equal to $-\frac 12 \alpha ^2$.

\begin{remark}\rm
A Lie group of Segre type $\{21 \}$ with $\alpha=0$ or $\alpha=\beta$ is
locally symmetric if and only if $\beta=0$ (see \cite{Calvaruso07}). This shows that previous
examples are not locally symmetric.
\end{remark}

\subsubsection{\underline{\emph{Segre type $\{3\}$}}.}

Assume $L$ has a triple root of its minimal polynomial. Then
\[
L=\left(\begin{array}{ccc}
\alpha &\frac{1}{\sqrt{2}} &\frac{1}{\sqrt{2}}\\
    \frac{1}{\sqrt{2}}   &\alpha&0\\
    -\frac{1}{\sqrt{2}}   & 0&\alpha
\end{array}\right)
\]
with respect to an orthonormal basis $\{e_1,e_2,e_3\}$ of signature
$(++-)$, and the corresponding Lie algebra is given by
\[
(\mathfrak{g}_{III}):\left\{\!\begin{array}{rcl}
[e_1,e_2]\!&\!=\!&\!-\frac{1}{\sqrt{2}}e_1-\alpha e_3, \qquad \quad
[e_1,e_3]\!\!=\!\!-\frac{1}{\sqrt{2}} e_1-\alpha e_2,\\
\noalign{\smallskip} [e_2,e_3]\!&\!=\!&\!\alpha
e_1+\frac{1}{\sqrt{2}}e_2 -\frac{1}{\sqrt{2}}e_3.
\end{array}\right.
\]
Hence the non-zero components of the curvature tensor (up to
symmetries) are
\[
\begin{array}{lll}
R_{1221}=\frac{1}{4} \left(\alpha ^2+4\right), &\quad
R_{1331}=1-\frac{1}{4}\alpha ^2, &\quad R_{2323}=\frac{1}{4}\alpha
^2,
\\[0.05in]
R_{1231}=1, &\quad R_{1223}=R_{1323}=\frac{1}{\sqrt{2}}\alpha. &
\end{array}
\]
The Ricci operator, expressed in terms of the basis $\{ e_1,e_2,e_3\}$,
becomes
\[
\hat{Ric}=\left(\begin{array}{ccc}
-\frac{1}{2}\alpha^2 &-\frac{1}{\sqrt{2}}\alpha&-\frac{1}{\sqrt{2}}\alpha\\[0.05in]
-\frac{1}{\sqrt{2}}\alpha&-\frac{1}{2}{(\alpha^2+2)}&-1\\[0.05in]
\frac{1}{\sqrt{2}}\alpha &1 &1-\frac{1}{2}\alpha^2
\end{array}\right),
\]
with a single eigenvalue $-\frac{1}{2}\alpha^2$.
If $\alpha\neq 0$, then $\hat{Ric}$ is of Segre type $\{3\}$,
while $\hat{Ric}$  is two-step nilpotent if $\alpha=0$.

For a vector $X=\sum X_ie_i$, the Lie derivative {has the following expression}
\[
(\mathcal{L}_Xg)= \frac{1}{\sqrt{2}}\left(
\begin{array}{ccc}
-2(X_2+X_3) & {X_1} & {X_1}\\
{X_1} & 2 X_3 & X_3-X_2\\
{X_1} & X_3-X_2& -2 X_2
\end{array}\right).
\]
Thus the Ricci soliton condition \eqref{soliton} on $\mathfrak{g}_{III}$ gives rise to the following system:
\begin{equation}\label{IIIeq}
\left\{
\begin{array}{l}
\frac{\alpha^2}{2}+\sqrt{2} X_2 +\sqrt{2} X_3=-\lambda{,}
\\
\noalign{\smallskip}
\frac{\alpha^2}{2}-\sqrt{2} X_3 + 1 = -\lambda{,}\\
\noalign{\smallskip}
\frac{\alpha^2}{2}-\sqrt{2} X_2-1= -\lambda{,}
\\
\noalign{\smallskip}
\frac{1}{\sqrt{2}}(X_1-\alpha)=0{,}
\\
\noalign{\smallskip}
X_2-X_3+\sqrt{2}=0{.}
\end{array}
\right.
\end{equation}

If we subtract half of the second and third equations to the first one in \eqref{IIIeq},
we see that $X_2=-X_3$ and therefore
$\lambda=-\frac{1}{2}\alpha^2$. Moreover from the forth equation $X_1=\alpha$. Hence \emph{any Segre type $\{3\}$
unimodular Lie group is a homogeneous Ricci soliton}
for
\begin{equation}\label{XIII}
X=\alpha\, e_1 -\frac{1}{\sqrt{2}}\, e_2 +\frac{1}{\sqrt{2}}\, e_3.
\end{equation}

\begin{remark}\rm
A vector field $X$ defining a homogeneous Ricci soliton on $\mathfrak{g}_{III}$ satisfies $\langle
X,X\rangle =\alpha^2$ and thus it is either spacelike or null. Correspondingly, the homogeneous Ricci soliton is either
expanding or steady.
Note also that Segre type $\{3\}$ unimodular Lie groups are never symmetric (see also \cite{Calvaruso07}).
\end{remark}

Previous analysis proves {\em iv)} of Theorem \ref{th:1}.
Lie groups having  unimodular Lie algebras compatible with the Ricci soliton equation \eqref{soliton},
and listed in Theorem \ref{th:1}, can be deduced from \cite{Rahmani} (se also \cite{calvaruso}).
The results we proved are summarized in the following

\begin{theorem}\label{unim}
The following are all  non-trivial homogeneous Lorentzian Ricci solitons realized as unimodular Lorentzian
Lie groups $G$:
\begin{itemize}
\item[a)] $G=E(1,1)$, with Lie algebra as in Theorem \ref{th:1}--{\em (iv.1)}, $\alpha=0\neq \beta$.
The homogeneous Ricci soliton is steady and defined by a spacelike vector field~\eqref{XIIa}.
\item[b)] $G=O(1,2)$ or $SL(2,\mathbb R)$, with Lie algebra as in Theorem \ref{th:1}--{\em (iv.1)}, $\alpha=\beta\neq 0$.
The homogeneous Ricci soliton is expanding and defined by a spacelike vector field~\eqref{XIIb}.
\item[c)] $G=O(1,2)$ or $SL(2,\mathbb R)$, with Lie algebra as in Theorem \ref{th:1}--{\em (iv.2)}, $\alpha\neq 0$.
The homogeneous Ricci soliton is expanding and defined by a spacelike vector field \eqref{XIII}.
\item[d)] $G=E(1,1)$, with Lie algebra as in Theorem \ref{th:1}--{\em (iv.2)}, $\alpha=0$.
The homogeneous Ricci soliton is steady and defined by a null vector field \eqref{XIII}.
\end{itemize}
\end{theorem}

\begin{remark}\rm
Ricci solitons listed in Theorem \ref{unim} are locally conformally flat if and only if they
correspond to $G=E(1,1)$, which is not locally symmetric.
\end{remark}

\subsection{Non-unimodular Lie groups}

Following \cite{CP-1}, non-unimodular Lorentzian Lie algebras of non-constant
sectional curvature are given, with respect to a suitable basis $\{
e_1,e_2,e_3\}$, by
\begin{equation}\label{IV}
(\mathfrak{g}_{IV}):\qquad \qquad [e_1,e_2] = 0, \quad [e_1,e_3] =
\alpha e_1 + \beta e_2, \quad [e_2,e_3] = \gamma e_1 + \delta e_2,
\end{equation}
where $\alpha+\delta\neq 0$ and one of the following holds:
\begin{enumerate}
\item[IV.1] $\{ e_1,e_2,e_3\}$ is orthonormal with $\langle e_1,e_1\rangle=-\langle e_2, e_2\rangle=-\langle e_3,e_3\rangle=-1$
and the structure constants satisfy $\alpha\gamma-\beta\delta=0$.
\item[IV.2] $\{ e_1,e_2,e_3\}$ is orthonormal with $\langle e_1,e_1\rangle=\langle e_2, e_2\rangle=-\langle e_3,e_3\rangle=1$
and the structure constants satisfy $\alpha\gamma+\beta\delta=0$.
\item[IV.3] $\{ e_1,e_2,e_3\}$ is a pseudo-orthonormal basis with
\[
\langle\,\cdot\,,\,\cdot\,\rangle=\left(\begin{array}{ccc}
1&0&0\\
0&0&-1\\
0&-1&0
\end{array}\right)
\]
and the structure constants satisfy $\alpha\gamma=0$.
\end{enumerate}

We analyze the three cases separately.

\subsubsection{\underline{\emph{Type IV.1}}.}

The non-zero components of the curvature tensor are given by
\[
\begin{array}{l}
R_{1212}=\frac{1}{4} \left(\beta ^2+\gamma^2 +4\alpha\delta -2\beta\gamma\right), \\[0.05in]
R_{1313}=\frac{1}{4} \left(4\alpha^2 -3\beta^2 +\gamma^2 +2\beta\gamma\right), \\[0.05in]
R_{2332}=\frac{1}{4} \left(\beta^2 -3\gamma^2 +4\delta^2
+2\beta\gamma\right).
\end{array}
\]
Hence the Ricci operator is diagonalizable with eigenvalues
\[
\begin{array}{rcl}
\lambda_1&=&
\frac{1}{2}(\beta^2-\gamma^2-2\alpha(\alpha+\delta)),\\
\noalign{\smallskip}
\lambda_2&=&
\frac{1}{2}(\gamma^2-\beta^2-2\delta(\alpha+\delta)),\\
\noalign{\smallskip}
\lambda_3&=& \frac{1}{2}((\beta-\gamma)^2-2(\alpha^2+\delta^2)).
\end{array}
\]
The Lie derivative of the metric for an arbitrary vector $X=\sum
X_ie_i$ is given by
\[
(\mathcal{L}_Xg)= \left(
\begin{array}{ccc}
-2\alpha X_3 & X_3(\beta-\gamma) & X_1\alpha+X_2\gamma\\
X_3(\beta-\gamma) & 2 X_3\delta & -X_1\beta-X_2\delta\\
X_1\alpha+X_2\gamma & -X_1\beta-X_2\delta & 0
\end{array}\right){,}
\]
and thus necessary and sufficient conditions for the existence of a
homogeneous Ricci soliton \eqref{soliton} on
$\mathfrak{g}_{IV.1}$ are given by
\begin{equation}\label{IV.1eq}
\left\{\!\begin{array}{l}
\beta^2-\gamma^2-2\alpha(\alpha+\delta)+4X_3\alpha=2\lambda{,}
\\
\noalign{\smallskip}
\gamma^2-\beta^2-2\delta(\alpha+\delta)+4X_3\delta=2\lambda{,}\\
\noalign{\smallskip}
(\beta-\gamma)^2-2(\alpha^2+\delta^2)=2\lambda{,}
\\
\noalign{\smallskip}
X_1\alpha+X_2\gamma=0{,}\\
\noalign{\smallskip}
X_1\beta+X_2\delta=0{,}
X_3(\beta-\gamma)=0.
\end{array}
\right.
\end{equation}

If $X_3=0$, then the first three equations in (\ref{IV.1eq}) imply
$\lambda_1=\lambda_2=\lambda_3$.
On the other hand, if $X_3 \neq 0$, then the last equation in (\ref{IV.1eq}) gives $\beta=\gamma$.
Since $\alpha\gamma-\beta\delta=0$, from $\alpha+\delta\neq 0$ and (\ref{IV.1eq}) we obtain
$\alpha=\delta$ and hence $\lambda_1=\lambda_2=\lambda_3$.
Thus,  all solutions of (\ref{IV.1eq}) are Einstein, and hence of constant sectional curvature.

\subsubsection{\underline{\emph{Type IV.2}}.}
Assume the non-unimodular Lie algebra $\mathfrak{g}_{IV}$ has a basis as given in IV.2. Then, a
straightforward calculation shows that the non-zero components of
the curvature tensor are given by
\[
\begin{array}{l}
R_{1212}=\alpha \delta -\frac{1}{4}(\beta +\gamma )^2, \\[0.05in]
R_{1331}=\frac{1}{4} \left(4\alpha^2 +3\beta^2 -\gamma^2 +2\beta\gamma\right), \\[0.05in]
R_{2323}=\frac{1}{4} \left(\beta^2 -3\gamma ^2 -4\delta^2
-2\beta\gamma\right).
\end{array}
\]
Therefore, the Ricci operator is diagonalizable with eigenvalues
\[
\begin{array}{rcl}
\lambda_1&=&\frac{1}{2}(\beta^2-\gamma^2+ 2\alpha(\alpha+\delta)),\\
\noalign{\smallskip}
\lambda_2&=&\frac{1}{2}(\gamma^2-\beta^2+2\delta(\alpha+\delta)),\\
\noalign{\smallskip}
\lambda_3&=&\frac{1}{2}((\beta+\gamma)^2+2(\alpha^2+\delta^2)).
\end{array}
\]
A straightforward calculation from (\ref{IV}), using the fact that
the structure constants satisfy $\alpha\gamma+\beta\delta=0$ and
$\alpha+\delta\neq 0$, shows that the Lie derivative of the metric
with respect to a vector $X=\sum X_i e_i$ is given by
\[
(\mathcal{L}_Xg)=\left(\begin{array}{ccc}
2 X_3\alpha & X_3(\beta+\gamma) & -X_1\alpha-X_2\gamma\\
X_3(\beta+\gamma) & 2 X_3\delta & -X_1\beta-X_2\delta\\
-X_1\alpha-X_2\gamma & -X_1\beta-X_2\delta & 0
\end{array}\right).
\]
Then  a Ricci soliton must satisfy
\begin{equation}\label{IV.2eq}
\left\{\!
\begin{array}{ll}
\beta^2-\gamma^2+ 2\alpha(\alpha+\delta)+4X_3\alpha=2\lambda{,}
\\
\noalign{\smallskip}
\gamma^2-\beta^2+2\delta(\alpha+\delta)+4X_3\delta=2\lambda{,}\\
\noalign{\smallskip}
(\beta+\gamma)^2+2(\alpha^2+\delta^2)=2\lambda{,}
\\
\noalign{\smallskip}
X_1\alpha+X_2\gamma=0{,}\\
\noalign{\smallskip}
X_1\beta+X_2\delta=0{,}
\\
\noalign{\smallskip}
X_3(\beta+\gamma)=0.
\end{array}
\right.
\end{equation}

A similar analysis to that developed for type IV.1 shows that
homogeneous Ricci solitons of type IV.2 necessarily are of constant
sectional curvature.

\subsubsection{\underline{\emph{Type IV.3}}.}
Let now $\mathfrak{g}_{IV}$ admit a pseudo-orthonormal basis as in IV.3. We then consider the orthonormal basis
\begin{center}
$\tilde e_1:= e_1,\quad \tilde
e_2:=\frac{1}{\sqrt{2}}(e_2-e_3),\quad \tilde
e_3:=\frac{1}{\sqrt{2}}(e_2+e_3),$
\end{center}
with signature $(++-)$. Then the non-zero components of the
curvature tensor are given by
\[
\begin{array}{ll}
R_{1212}=\frac{1}{4}(2\alpha\delta-2\alpha^2-\gamma(2\beta+\gamma)),\,\,
& R_{1213}=\frac{1}{2}(\alpha^2+\beta\gamma-\alpha\delta),\,\,\\
\noalign{\medskip}
R_{1313}=\frac{1}{4}(2\alpha\delta-2\alpha^2+\gamma(\gamma-2\beta)),\,\,
& R_{2323}=-\frac{3}{4}\gamma^2.
\end{array}
\]

The Ricci operator in the new basis $\{\tilde e_1,\tilde e_2,\tilde
e_3\}$ becomes
\[
\hat{Ric}= \left(
\begin{array}{ccc}
-\frac{1}{2}\gamma^2 &0 &0\\
0&\frac{1}{2}\left(\alpha(\delta-\alpha)+\gamma(\gamma-\beta)\right)
&\frac{1}{2}\left(\alpha(\alpha-\delta)+\beta\gamma\right)\\
0& -\frac{1}{2}\left(\alpha(\alpha-\delta)+\beta\gamma\right)
&\frac{1}{2}\left({
\alpha(\alpha-\delta)+\gamma(\beta+\gamma)}\right)
\end{array}
\right),
\]
which has eigenvalues $\lambda_1=-\frac{1}{2}\gamma^2$ and
$\lambda_2=\lambda_3=\frac{1}{2}\gamma^2$.
Thus, $\hat{Ric}$ is of Segre type $\{21\}$.

For an arbitrary vector $X=\sum X_i \tilde{e}_i$, the Lie derivative
of the metric becomes
\[
(\mathcal{L}_Xg)=\left(\begin{array}{ccc}
\sqrt{2}\alpha(X_3-X_2) & a_{12}
& a_{13}\\
a_{12}
&X_1\beta+\sqrt{2}X_3\delta
& -X_1\beta-\frac{1}{\sqrt{2}}\delta (X_2+X_3)\\
a_{13}
& -X_1\beta-\frac{1}{\sqrt{2}}\delta (X_2+X_3) &
X_1\beta+\sqrt{2}X_2\delta
\end{array}\right),
\]
where
\[
\begin{array}{l}
a_{12}=\frac{1}{2}\left( X_3(\beta+2\gamma)+\sqrt{2}X_1\alpha-X_2\beta  \right),
\\
\noalign{\smallskip}
a_{13}=\frac{1}{2}\left( X_2(\beta-2\gamma)-\sqrt{2}X_1\alpha-X_3\beta  \right).
\end{array}
\]
Now, since $\alpha\gamma=0$, we consider the possibilities
$\alpha=0$ and $\gamma=0$ separately. Note that if $\alpha=\gamma=0$, then the
metric $\mathfrak{g}_{IV}$ is flat.

Assume first that $\alpha=0\neq \gamma$. Then,
(\ref{soliton}) holds if and only if
\begin{equation}\label{IVeq}
\left\{
\begin{array}{ll}
\gamma^2=-2\lambda{,}
\\
\noalign{\smallskip}
\gamma^2-\beta\gamma +2X_1\beta+2\sqrt{2}X_3\delta=2\lambda{,}\\
\noalign{\smallskip}
\gamma^2+\beta\gamma-2X_1\beta-2\sqrt{2}X_2\delta=2\lambda{,}
\\
\noalign{\smallskip}
-2\gamma X_3+X_2\beta-X_3\beta =0{,}\\
\noalign{\smallskip}
2\gamma X_2-X_2\beta+X_3\beta=0{,}
\\
\noalign{\smallskip}
\beta\gamma - 2X_1\beta-\sqrt{2}X_2\delta-\sqrt{2}X_3\delta=0{,}
\end{array}
\right.
\end{equation}
for a vector $X=\sum X_i\tilde{e}_i$. From the forth and fifth equation in \eqref{IVeq} we get
$X_2=X_3$, which  implies, using the second and third equations, that \eqref{IVeq} admits no solutions.

Now assume now $\alpha\neq 0= \gamma$. Then,
(\ref{soliton}) reduces to the following system of equations
\begin{equation}\label{IVeqb}
\left\{
\begin{array}{ll}
{2\sqrt{2}\alpha(X_2-X_3)}=-2\lambda{,}
\\
\noalign{\smallskip}
\alpha\delta-\alpha^2+2X_1\beta+2\sqrt{2}X_3\delta=2\lambda{,}\\
\noalign{\smallskip}
\alpha^2-\alpha\delta-2X_1\beta-2\sqrt{2}X_2\delta=2\lambda{,}
\\
\noalign{\smallskip}
\sqrt{2}X_1\alpha - X_2\beta +X_3\beta=0{,}\\
\noalign{\smallskip}
\alpha^2-\alpha\delta - 2X_1\beta=\sqrt{2}\delta(X_2+X_3){,}
\end{array}
\right.
\end{equation}
for a vector $X=\sum X_i\tilde{e}_i$.
We subtract the third equation to the second one and conclude, using the first equation, that either
$2\alpha=\delta$ or $X_2=X_3$. We analyze both cases separately.

Set first $\alpha=\frac{1}{2}\delta\neq 0$.
Then there exists homogeneous Ricci
solitons for
\begin{equation}\label{XIVa}
X= -\frac{2\beta\lambda}{\delta^2}\,\, \tilde e_1
-\frac{\delta^4+8(\delta^2-2\beta^2)\lambda}{8\sqrt{2}\delta^3}\,\, \tilde e_2
-\frac{\delta^4-8(\delta^2+2\beta^2)\lambda}{8\sqrt{2}\delta^3}\,\, \tilde e_3.
\end{equation}
Note that the  corresponding solitons may be expanding, steady or shrinking
depending on the value of $\lambda$, which can be chosen with
absolute freedom.

Set now $X_2=X_3$. Then necessarily $\lambda=0$ and $X_1=0$.
The remaining equation $\alpha^2-\alpha\delta -2\sqrt{2}X_2\delta=0$  in \eqref{IVeqb} gives
rise to homogeneous  steady Ricci solitons for
\begin{equation}\label{XIVb}
X=\frac{\alpha^2-\alpha\delta}{2\sqrt{2}\delta}\, \left(\tilde e_2+\tilde e_3\right).
\end{equation}

In all cases above the Ricci operator is two-step nilpotent and the
metric is non-symmetric whenever $\alpha\delta(\alpha-\delta)\neq 0$.
Furthermore, for the particular choice  $\delta=0 \neq \alpha$
the resulting metric is symmetric but not of constant curvature \cite{Calvaruso07}.

The results of this subsection prove case {\em v)} of Theorem \ref{th:1} and are summarized in the following

\begin{theorem}\label{nonunim}
A non-unimodular Lie group $G$ equipped with a left-invariant Lorentz\-ian metric
is a non-trivial homogeneous Ricci soliton if and only if its non-unimodu\-lar Lie algebra
$\mathfrak{g}_{IV}$ satisfies $\alpha \neq 0 =\gamma$.

Steady Ricci solitons, defined by null vector fields \eqref{XIVb},
exist for any choice of $\alpha\neq 0,\beta,\delta$.

In the special case $\delta=2\alpha$, there exist expanding, steady and shrinking Ricci solitons,
defined by vector fields \eqref{XIVa}, whose causal character depends on $\lambda$.
\end{theorem}

\begin{remark}\rm
Nontrivial Ricci solitons in Theorem \ref{nonunim} are locally conformally flat if and only if
$\gamma=\beta=0$, in which case they are not locally symmetric.
Therefore, Theorem \ref{nonunim} provides examples of complete locally conformally flat
expanding, steady and shrinking Ricci solitons.
\end{remark}

\section{Three-dimensional Walker manifolds}\label{s:walker}

We now consider three-dimensional Lorentzian manifolds
$(M,g)$ admitting a parallel null vector field $\mathcal{U}$.
We refer to \cite{walker-metrics,C-GR-VA}
and references therein for more information on Walker manifolds.
It has been shown by Walker \cite{W-50} that there exist adapted
coordinates $(t,x,y)$
where the Lorentzian metric tensor expresses as
\begin{equation}\label{eq:metric}
g=\left(\begin{array}{ccc}
0 & 0 & 1\\
0 & \varepsilon & 0\\
1 & 0 & f(x,y)
\end{array}\right),
\end{equation}
for some function $f(x,y)$, where $\varepsilon=\pm 1$ and the
parallel null vector field is ${\mathcal U}=\frac{\partial}{\partial t}$.
Then, the associated
Levi-Civita connection  is described by
\begin{equation}\label{eq:Levi-Civita}
\begin{array}{l}
\nabla_{\partial_x}\partial_y=\frac{1}{2}f_x\partial_t, \qquad
\nabla_{\partial_y}\partial_y=\frac{1}{2}f_y\partial_t
-\frac{1}{2\varepsilon}f_x\partial_x.
\end{array}
\end{equation}
As shown in \cite{C-GR-VA}, the Ricci tensor $Ric$ and the
Ricci operator $\hat{Ric}$ of
a metric (\ref{eq:metric}), expressed in the coordinate basis, take the form
\begin{equation}\label{eq:Ric}
Ric=-\frac{1}{2\varepsilon}f_{xx}\left(\begin{array}{ccc}
0 & 0 & 0\\
\noalign{\smallskip}
0 & 0 & 0\\
\noalign{\smallskip}
0 & 0 & 1
\end{array}\right), \qquad \hat{Ric}=-\frac{1}{2\varepsilon}f_{xx}\left(\begin{array}{ccc}
0 & 0 & 1\\
\noalign{\smallskip}
0 & 0 & 0\\
\noalign{\smallskip}
0 & 0 & 0
\end{array}\right).
\end{equation}
Hence, if $f_{xx}=0$, then the metric \eqref{eq:metric} is flat,
while for $f_{xx} \neq 0$ the Ricci operator $\hat{Ric}$ is two-step nilpotent.

Next, let $X=(A(t,x,y),B(t,x,y),C(t,x,y))$ be an arbitrary vector field on $M$.
A straightforward calculation
from (\ref{eq:Levi-Civita})
shows that the Lie derivative of the metric $(\mathcal{L}_Xg)$
expresses in the coordinate basis as follows:
\begin{equation}\label{eq:lied}
\mathcal{L}_Xg\!=\!\!\left(\!
\begin{array}{ccc}
2 C_t & \varepsilon B_t+C_x &A_t+C_y+fC_t \\
\noalign{\smallskip}
\varepsilon B_t+C_x & 2 \varepsilon B_x & A_x+\varepsilon B_y+fC_x\\
\noalign{\smallskip}
A_t+C_y+fC_t & A_x+\varepsilon B_y+fC_x & B f_x +  C f_y +2(A_y+fC_y)
\end{array}
\!\right)\!.
\end{equation}

Now, from (\ref{eq:Ric}) and (\ref{eq:lied}) we obtain the following necessary and
sufficient
conditions for a strict Walker metric \eqref{eq:metric} to be a Ricci soliton:
\begin{equation}\label{walker-eq}
\left\{
\begin{array}{l}
2 C_t=0,\\
\noalign{\smallskip}
C_x+\varepsilon B_t=0,\\
\noalign{\smallskip}
C_y+A_t+fC_t=\lambda,\\
\noalign{\smallskip}
2B_x=\lambda,\\
\noalign{\smallskip}
A_x+\varepsilon B_y + fC_x=0,\\
\noalign{\smallskip}
2A_y+f_xB+f_yC-\frac{1}{2\varepsilon}f_{xx}=f(\lambda-2C_y).
\end{array}
\right.
\end{equation}

The first equation in \eqref{walker-eq} gives $C=C(x,y)$
and simplifies the third one. Since $C$ does not depend on $t$,
we can  easily integrate the second and third equations in \eqref{walker-eq}
and get
\[
A=(\lambda-C_y)t+ G(x,y), \qquad B={-\frac{C_x t}{\varepsilon}}
+H(x,y).
\]
Therefore, the fourth equation in \eqref{walker-eq} now gives
\begin{equation}\label{IIw}
-2tC_{xx}  +2\varepsilon H_x = \varepsilon\lambda.
\end{equation}
Since \eqref{IIw} must hold for any value of $t$,
it implies at once $C_{xx}=0$ and $2H_x =\lambda$.
By integration we then have $C=u(y)x+v(y)$ and $H=\frac 12 \lambda x + w(y)$.
Then the fifth equation in \eqref{walker-eq} becomes
\[
f u(y) - 2 t u'(y)+\varepsilon w'(y)+G_x=0
\]
from where it follows that $u(y)$ is constant: $u(y)=\alpha$.
Then, system \eqref{walker-eq} now reduces to
\begin{equation}\label{walker-eq2}
\left\{
\begin{array}{l}
\alpha f +\varepsilon w'(y)+G_x=0,\\
\noalign{\smallskip}
2f v'(y)-\lambda f- 2tv''(y)
+f_y(\alpha x+v(y))+2G_y\\
\noalign{\smallskip} \phantom{2f v'(y)-\lambda f- 2tv''(y)}
+f_x(w(y)+\frac{\lambda}{2}x-\varepsilon\alpha t)
=\frac{1}{2\varepsilon}f_{xx}.
\end{array}
\right.
\end{equation}

The second equation in \eqref{walker-eq2} holds for any value of
$t$. Hence ${\frac{\alpha f_x}{\varepsilon}}=-2v''(y)$
and thus $\alpha f_{xx}=0$. Since $f_{xx}=0$ if and only if the
Walker metric is flat, we assume $\alpha=0$ and \eqref{walker-eq2}
reduces to
\begin{equation}\label{walker-eq3}
\left\{
\begin{array}{l}
\varepsilon w'(y)+G_x=0,\\
\noalign{\smallskip}
2f v'(y)-\lambda f- 2tv''(y)
+f_yv(y)+2G_y
+f_x(w(y)+\frac{\lambda}{2}x)
=\frac{1}{2\varepsilon}f_{xx}.
\end{array}
\right.
\end{equation}
The first equation in \eqref{walker-eq3} gives
$G(x,y)=-\varepsilon x w'(y)+\mu(y)$. Then, since the second equation in \eqref{walker-eq3}
must be independent of $t$, one gets $v''(y)=0$ and hence $v(t)=\beta y + \gamma$.

Finally we conclude that there exist nontrivial Ricci solitons given by Walker
metrics \eqref{eq:metric} if and only if the vector field $X$ takes the form
$$
X(t,x,y)\!=\! \left(
t(\lambda-\beta)-\varepsilon x w'(y)+\mu(y),
\frac{1}{2} \lambda x + w(y),
\beta y+\gamma
\right),
$$
for some real constants $\beta,\gamma$ and smooth functions $w(y)$
and  $\mu(y)$, satisfying the partial differential equation
\begin{equation}\label{walker-eq4}
2\beta f -\lambda f+2\mu'(y)-2\varepsilon xw''(y)
+f_y(\beta y +\gamma)
+f_x(\frac{\lambda}{2}x+w(y))
=\frac{1}{2\varepsilon}f_{xx}.
\end{equation}

One can not expect the partial differential equation \eqref{walker-eq4}
to admit solutions in general.
We now turn our attention to the special case when the Walker metric
is locally symmetric.
Locally symmetric Walker metrics \eqref{eq:metric} are characterized
by the fact that their defining function $f(x,y)$ is given by (see \cite{walker-metrics, C-GR-VA})
\begin{equation}\label{lsWalker}
f(x,y)= x^2 \kappa+x P(y)+Q(y),
\end{equation}
for arbitrary functions $P$ and $Q$, and constant $\kappa$ which vanishes if
and only if the metric is flat.
When $f$ satisfies \eqref{lsWalker}, a straightforward calculation
shows that equation \eqref{walker-eq4} becomes
\begin{equation}\label{wls}
\begin{array}{l}
0=2x^2\beta\kappa +x\left( 2\beta P(y)-\frac{1}{2}\lambda P(y)+2\kappa w(y)
+(\beta y+\gamma)P'(y)-2\varepsilon w''(y)\right)\\
\noalign{\smallskip}
\phantom{0=}
-\frac{\kappa}{\varepsilon}+2\beta Q(y)-\lambda Q(y)+P(y)w(y)
+(\beta y+\gamma)Q'(y)+2\mu'(y).
\end{array}
\end{equation}
But \eqref{wls} must hold for all values of $x$.
Therefore, it gives $\beta=0$ (excluding the flat case $\kappa=0$)
and reduces to the system
\begin{equation}\label{wfinal}
\left\{
\begin{array}{l}
2\varepsilon w''(y)-2\kappa w(y)=\gamma P'(y)-\frac{1}{2}\lambda P(y),\\
\noalign{\smallskip}
2\mu'(y)=\frac{\kappa}{\varepsilon}-P(y)w(y)+\lambda Q(y)-\gamma
Q'(y).
\end{array}
\right.
\end{equation}

The second equation in \eqref{wfinal}, by direct integration,
permits to express $\mu(y)$ in terms of $w(y)$ and $P(y),Q(y)$.
The first equation in \eqref{wfinal} is a second order linear ordinary
differential equation for $w(y)$, with constant coefficients,
determined by the smooth function $\gamma P'(y)-\frac{\lambda}{2} P(y)$.
A standard theorem ensures the existence of solutions for such an equation.
Therefore, we proved the following

\begin{theorem}\label{wsolitons}
Any three-dimensional symmetric Walker metric \eqref{eq:metric} is a  Ricci soliton,
which can be expanding, steady or shrinking and is defined by
vector fields
\begin{equation}\label{eq:ws-1}
X(t,x,y)\!=\! \left(
\lambda t-\varepsilon xw'(y)+\mu(y),
\frac{1}{2}\lambda x + w(y),
\gamma
\right),
\end{equation}
where $\lambda$ and $\gamma$ are real constants and the functions
$w(y)$ and $\mu(y)$ are arbitrary solutions of \eqref{wfinal}.
In general, the causal character of $X$ may vary with the point.
\end{theorem}

%

\end{document}